\documentclass[1p,numbers]{elsarticle}
\usepackage{amsmath,amsfonts,enumerate,amssymb,amscd,scrpage2}
\usepackage[numbers,square]{natbib}
\usepackage[thmmarks,amsmath]{ntheorem}
\usepackage[arrow, matrix, curve]{xy}

\newtheorem{theorem}{Theorem}[section]
\newtheorem{definition}[theorem]{Definition}
\newtheorem{lemma}[theorem]{Lemma}

\newtheorem{proposition}[theorem]{Proposition}

\theoremstyle{nonumberplain}
\theoremsymbol{\ensuremath{\square}}
\theoremheaderfont{\scshape}
\theoremseparator{:}
\theorembodyfont{\normalfont}
\newtheorem{proof}{Proof}

\begin{document}
\title{Hedlund Metrics and the Stable Norm}
\author{Madeleine Jotz}
\ead{madeleine.jotz@epfl.ch}
\address{Section de Math{\'e}matiques, Ecole Polytechnique F{\'e}d{\'e}rale de Lausanne,
 1015 Lausanne, Switzerland.}

\begin{abstract}
The real homology of a compact Riemannian manifold $M$ is
  naturally endowed with the stable norm. The stable norm on $H_1(M,\mathbb{R})$ arises
  from the Riemannian length functional by homogenization. It is difficult and
  interesting to decide which norms on the finite-dimensional vector space
  $H_1(M,\mathbb{R})$ are stable norms of a Riemannian metric on $M$. If the dimension of
  $M$ is at least three, I. Babenko and F. Balacheff proved in \cite{baba}
  that every polyhedral norm ball in $H_1(M,\mathbb{R})$, whose vertices are rational
  with respect to the lattice of integer classes in $H_1(M,\mathbb{R})$, is the stable
  norm ball of a Riemannian metric on $M$. This metric can even be chosen to
  be conformally equivalent to any given metric. In \cite{baba}, the stable
  norm induced by  the constructed metric is computed by  comparing the metric with a polyhedral one.  Here we
  present an alternative construction for the  metric, which remains in the
  geometric framework of smooth Riemannian metrics.
\end{abstract}

\begin{keyword}
 Riemannian metrics, stable norm, polytopes.
\MSC[2008] 53C22, 53C38, 58A10, 58F17, 53B21
\end{keyword}

\maketitle
\section{Introduction}
On every compact Riemannian manifold $M$ the real homology vector spaces
$H_m(M;\mathbb{R})$ are endowed with a natural norm $\|\cdot\|_s$, called \emph{stable
  norm}. This concept appeared  for the first time in \citet{federer} and was
named \emph{stable norm} in \citet{gromov}. 
The stable norm on $H_1(M;\mathbb{R})$ arises directly from the Riemannian metric on the
manifold $M$.
The following equality for an integral class $v\in H_1(M;\mathbb{R})$ (see \cite{gromov})
\begin{equation*}
\|v\|_s:=\inf\{n^{-1}L(\gamma)| \gamma \text{ is a closed curve representing }nv,\,n \in \mathbb{N}\} 
\end{equation*} allows a description of this object that is geometrically very
intuitive: the stable norm describes the geometry of the Abelian covering
$\bar{M}$ of $M$ from a point of view from which fundamental domains look
arbitrarily small.
Knowing the unit ball of this norm, one can decide on existence and properties
of some of the minimal geodesics relative to the Riemannian Abelian covering of the
manifold; these are curves in $M$ whose lifts to the Riemannian
Abelian  covering minimize arc length between each two of their points.
Bangert has presented  in \citep{ba} a Riemannian metric on the $3$-torus
$\mathbb{T}^3$, such that the unit ball of the induced stable norm on $H_1(\mathbb{T}^3;\mathbb{R})\simeq \mathbb{R}^3$ is a symmetric
octahedron. Furthermore, Babenko and Balacheff have shown in \citep{baba} that,
given a compact Riemannian manifold $(M,\rho)$ of dimension greater than $2$,
for every centrally symmetric and convex polytope in $H_1(M;\mathbb{R})$ with nonempty interior, such that the directions of its vertices
are rational, there is a Riemannian metric on
$M$ that is conformal to $\rho$ and induces the given polytope as unit ball
of the stable norm.  
Here we propose an alternative Riemannian metric, satisfying the same
conditions. Our construction is a generalization of the Hedlund metric in
\citet{ba}. The idea, that can be already found in the original paper of
\citet{hedlund} and is also used in \cite{baba}, is to construct a metric that is ``small'' in tubular
neighborhoods of disjoint closed curves representing the vertices of the
polytope, and much ``bigger'' everywhere else. 
The convexity properties of the
polytope play a decisive role in our computation of the stable norm induced
by the Hedlund metric.   

Bangert and Hedlund use such metrics
in order to illustrate their results on minimal geodesics. Here we
focuse only on the proof of the theorem of \citet{baba}. In fact, if we wanted
to show results on minimal geodesics, we would need to specify the definition
of the Hedlund metric we give here. A discussion of the minimal geodesics for
such metrics (with additional assumptions) was made in \citet{jotz}.

\paragraph{Outline of the paper:} in the next section the construction of
tubular neighborhoods of curves will be recalled. There a lemma on existence
of representatives for cohomology classes with ``good''
properties on the tubular neighborhood will be stated. In the following section
the construction of the Riemannian metric will be given and the formula for
the corresponding stable norm will be computed.

\paragraph{Notations:} In the following $M$ will denote a compact smooth manifold with $\dim M\geq 3$ and $\rho$  a Riemannian metric on $M$. Let $\bar M$
denote the Abelian covering of $M$. More precisely $\bar M$ is the subcovering
of the universal covering whose group of deck transformations is the set
$H_1(M;\mathbb{Z})_{\mathbb{R}}$ of integer classes in
$H_1(M;\mathbb{R})$. We denote by $p:\bar M\to M$
the covering map and by  $\bar \rho:= p^*\rho$ the pull-back metric.
If $h:\pi_1(M)\to H_1(M;\mathbb{Z})$ denotes the Hurewicz homomorphism
(\citep[see][]{leetop}) and $T$  the torsion subgroup of
$H_1(M;\mathbb{Z})$, then the Abelian covering
can be described as the quotient manifold of the action of the normal
subgroup $h^{-1}(T)\subseteq \pi_1(M)$ of the fundamental group on the universal
cover $\tilde M$  of $M$. Hence the operation 
\begin{align*}
\begin{array}{lccc}
\Phi:&H_1(M;\mathbb{Z})_{\mathbb{R}}\times\bar{M}&\to&\bar{M}\\
&(v,m)&\mapsto&\Phi(v,m)=:m+v
\end{array}
\end{align*}
of $H_1(M;\mathbb{Z})_{\mathbb{R}}$ on $\bar{M}$ is abelian and
torsionfree (that is why we choose to use this $+$-notation).\\
The de Rham cohomology
vector space $H^1_{\bf dR}(M)$ is isomorphic to the dual of $H_1(M;\mathbb{R})$
\citep[de Rham theorem]{lee}. In the following, we will use this isomorphism
without mentioning it.

Given a Riemannian metric $g$ on $M$, we will write $g^*$ for its dual metric.
The space of $1$-forms on $M$ (respectively on $\bar{M}$) will be denoted by
$\Omega^1(M)$ (respectively $\Omega^1(\bar{M})$).  We will denote by
$\|\cdot\|_x$ (or also simply $\|\cdot\|$) the norm on $T_xM$ induced by the
considered metric on $M$ (we will also use this notation for the norm on
$T_{\bar{x}}\bar{M}$, $\bar{x}\in\bar{M}$ induced from the corresponding
metric on $\bar{M}$). For a curve $\gamma:I\to M$, $L(\gamma)$ will be the
length  induced from the given metric on $M$ and for a curve $\bar\gamma:I\to
\bar{M}$, $\bar{L}(\bar\gamma)$ the length induced from the corresponding
periodic metric on $\bar{M}$.

Given a polytope $P$, we will call the set 
$\{\sum_{i=1}^k\alpha_iv_i\mid \alpha_i\geq 0\}$
the \emph{cone over the face $S$}  of the polytope, where $v_1,\dots,v_k$ are the
vertices of $P$ lying in this face (i.e. $S=\{\sum_{i=1}^k\alpha_iv_i\mid
\alpha_i\geq 0 \text{ and }  \sum_{i=1}^k\alpha_i=1\}$).

An integer class $v$ in $H_1(M;\mathbb{Z})_{\mathbb{R}}$ will be called \emph{indivisible} if the equation
$v=n\cdot v'$, $ n\in \mathbb{Z}$ and $v'\in H_1(M;\mathbb{Z})_{\mathbb{R}}$ yields $n=\pm 1$.

\paragraph{Acknowledgment:} I would like to thank Prof. Victor
Bangert who supervised my diploma thesis and gave me much advice for
this paper. I am also very grateful that he gave me the possibility to stay at the University of
Freiburg during a few months after my diploma.

I also thank the referees for many useful comments.  
\section{Tubular neighborhoods of curves, adapted one-forms}

\paragraph{Tubular neighborhoods and semi-geodesic coordinates.}
Let  $\gamma:[0,1]\to M$ be a regular simple closed curve. In the following, such a curve will be
called \emph{admissible}. We can write  $\gamma:\mathbb{S}^1\to
M$ and assume the curve $\gamma$ is parametrized proportionally to
arc length.
  
For $\varrho>0$ let $V_\varrho(\Gamma)$ denote the bundle of balls of radius
$\varrho$ in the normal bundle $\pi:N\Gamma\to\Gamma$ of the embedded
submanifold $\Gamma:=\gamma(\mathbb{S}^1)$ in
$M$. Analogously, if $I\subseteq\mathbb{S}^1$ is an interval, then
$V_\varrho(\gamma(I))=V_\varrho(\Gamma)\cap\pi^{-1}(\gamma(I))$. We choose $\varrho>0$
small enough such that the normal exponential map $E$ restricted to $V_\varrho(\Gamma)$
is a diffeomorphism onto an open neighborhood $U_\varrho(\Gamma)\subseteq M$
of $\Gamma$ (and similarly $U_\varrho(\gamma(I))=E(V_\varrho(\gamma(I)))$\,).
Such an open set  $U_\varrho(\Gamma)$ is called the \emph{tubular neighborhood} (of radius $\varrho$) of $\Gamma$.

Choose an orthogonal frame $(E_1,\ldots,E_m)$ on $U \subseteq M$ open, such
that for all $x=\gamma(t)$ in  $\Gamma\cap U$,
\[E_1\arrowvert_{x}=\dot{\gamma}(t)\] and, consequently,
$(E_2\arrowvert_{x},\ldots,E_m\arrowvert_{x})$ forms a basis for $N_x\Gamma$. Assume the open set
$U$ is such that $U_\varrho(\Gamma)\cap U=U_\varrho(\gamma(I))$ for an
open interval $I\subseteq \mathbb{S}^1$. The diffeomorphism \[\begin{array}{rlcc}\varphi:&U_\varrho(\gamma(I))&\to
&I\times B_\varrho^{m-1}\subseteq \mathbb{R}^m\\
&x&\mapsto&(s(x),\varphi_2(x),\dots,\varphi_m(x)),\end{array} \]
where $\varphi_j(x)$ and  $s(x)$ are such that
\begin{equation*}
E^{-1}(x)=\sum_{j=2}^m\varphi_j(x)\cdot
E_j\arrowvert_{\gamma(s(x))}\in \mathcal{V}_\varrho,
\end{equation*}
will be called a
  \emph{semigeodesic chart for $U_\varrho(\Gamma)$}. A particularity of this chart is that $\partial^\varphi_1\arrowvert_{x}=\dot{\gamma}(t)$ and, for $j=2,\ldots,m$,   
$\partial^\varphi_j\arrowvert_{x}=E_j\arrowvert_{x}$ holds for all $x=\gamma(t)\in
  \Gamma\cap U$ (note that $\Gamma\cap U=\gamma(I)$\,).

The map $s$ is defined globally on $U_\varrho(\Gamma)$ and we have the
identity
\begin{equation}
ds\arrowvert_{\gamma(t)}(\dot{\gamma}(t))=\frac{d}{dt}s\circ\gamma(t)=\frac{d}{dt}t=1\label{ds}
\end{equation}
for all $t$ in $\mathbb{S}^1$.

Let $\gamma_1,\dots,\gamma_N$ be  disjoint admissible loops and choose
$\varrho>0$ so that the construction above is possible  for all the curves
$\gamma_1,\dots,\gamma_N$ simultaneously. Choose furthermore $\varepsilon$
with $\varrho>\varepsilon>0$  such that  the tubular neighborhoods with radius
$\varepsilon$ of the curves
 are disjoint. Set $\Gamma_j=\gamma_j(\mathbb{S}^1)$, $\Gamma=\cup_{j=1}^N\Gamma_j$, and 
 $U_\varepsilon(\Gamma):=\cup_{j=1}^NU_\varepsilon(\Gamma_j)$. Then there exists a bump-function $\zeta$
on $M$ for the tubular neighborhoods, i.e., $\zeta$ is a smooth function such
that the following holds:

\begin{align}\label{stern} 
          \zeta(y)&=\left\{\begin{array}{ll}
             1&,\ y\in U_\varepsilon(\Gamma)\\
              0      & ,\ y \in M\setminus U_\varrho(\Gamma).\end{array}\right.   
\end{align}

\paragraph{``Good'' one-forms.}
Choose a connected fundamental domain $F_0$ for the action of $H_1(M;\mathbb{Z})_{\mathbb{R}}$ on $\bar{M}$. Denote by $\bar{\gamma}_j$ the lift
of $\gamma_j$ to $\bar{M}$ such that $\bar{\gamma}_j(0)\in F_0$ (note that
$\gamma_j$ is here considered as a smooth $1$-periodic curve
$\gamma_j:\mathbb{R}\to M$). Write
$\bar{\Gamma_i}=\bar{\gamma_i}(\mathbb{R})$ and $U_\varrho(\bar{\Gamma_i})$ the
corresponding lift to $\bar{M}$ of
$U_\varrho(\Gamma_i)$. Hence $U_\varrho(\bar{\Gamma_i})$ is the tubular
neighborhood of radius $\rho$ of $\bar{\Gamma_i}$. Thus the notion of a semigeodesic chart for $U_\varrho(\bar{\Gamma_i})$
makes also sense here, and $\bar{s}_i:U_\varrho(\bar{\Gamma_i})\to\mathbb{R}$ exists with 
$\bar{s}_i(\bar{\gamma_i}(t))=t$ for all $t\in\mathbb{R}$.\label{barsi}
Since the covering map $p:\bar M\to M$ is a local isometry,
\[\bar x\in \exp_{\bar{M}}(N_{\bar{\gamma_i}(t)}\bar{\Gamma_i})
\Leftrightarrow p(\bar x)
\in \exp_{M}(N_{p\circ\bar{\gamma}_i(t)}\Gamma_i)\] holds for all $\bar x\in
U_\varrho(\bar{\Gamma_i})$
  and 
\begin{equation}\label{blab}
(p^*ds_i)\arrowvert_{U_\varrho(\bar{\Gamma}_i)}=d\bar{s}_i. 
\end{equation}

Define  $L_i=\bar{\Gamma_i}+H_1(M;\mathbb{Z})_{\mathbb{R}}$ and
$U_\varrho(L_i)=U_\varrho(\bar{\Gamma_i})+H_1(M;\mathbb{Z})_{\mathbb{R}}$, as well as
$L=\cup_{j=1}^N L_j$ and $U_\varrho(L)=\cup_{j=1}^N
U_\varrho(L_j)$. Choose $\varepsilon$ with $0<\varepsilon<\varrho$ and define
$U_\varepsilon(\bar{\Gamma_i})$,
$U_\varepsilon(L_i)$ and $U_\varepsilon(L)$ as above. The connected
components of $L$ will be called \emph{lines} in the following.

In the following, a regular simple closed
curve will be called an \emph{admissible} curve.

\begin{proposition}\label{funk}
Let $v_1, \dots, v_N$ be indivisible integer classes in
$H_1(M;\mathbb{Z})_{\mathbb{R}}$, that span $H_1(M;\mathbb{R})$ as a real
vector space. Let $\gamma_1, \dots,
\gamma_N$ be disjoint admissible representatives of those
classes,  and
$U_\varepsilon(\Gamma_1),\dots,U_\varepsilon(\Gamma_N)$ disjoint
tubular neighborhoods of these curves.
Furthermore let $\lambda \in H^1_{dR}(M)$ be an arbitrary cohomology class. 
Then there exists a one-form $\omega$ representing $\lambda$ such that:
\begin{equation*}
\omega\arrowvert_{x}=\lambda(v_i)ds_i\arrowvert_{x} \text{ for } x \in
U_\varepsilon(\Gamma_i), \quad i=1,\dots, N.
\end{equation*}
\end{proposition}

\begin{proof}
 For $j=1,\dots,N$, the function $\bar{s}_j$ is defined on
  $U_\varrho(\bar{\Gamma}_j)$. Set $\bar{s}_j=0$ on
  $U_\varrho(\bar{\Gamma}_i)$ for $i\neq j$ and define:
\begin{align*}
s_\lambda:U_\varrho(L) & \to  \mathbb{R}\\
x=x_0+v_0 &\mapsto  \sum_{i=1}^N\lambda(v_i)\bar{s}_i(x_0)+\lambda(v_0). 
\end{align*}
Doing so, each element $U_\varrho(L_j)$ is written  $x=x_0+v_0$ with $x_0 \in
U_\varrho(\bar{\Gamma}_j)\cap F_0$  and
$v_0 \in H_1(M;\mathbb{Z})_{\mathbb{R}}$.
For $x \in
U_\varrho(\bar{\Gamma}_j)\cap F_0$ holds:
$s_\lambda(x)=\lambda(v_j)\bar{s}_j(x)$. Thus, with the definition of $s_\lambda$, for
$v=z\cdot v_j$ with $z \in \mathbb{Z}$:
\begin{align*}
s_\lambda(x+v)=\lambda(v_j)\bar{s}_j(x)+\lambda(v)
   =\lambda(v_j)\cdot(\bar{s}_j(x)+z)
\overset{\eqref{blab}}=\lambda(v_j)\cdot\bar{s}_j(x+v).
\end{align*}
This leads to
$s_\lambda\arrowvert_{U_\varrho(\bar{\Gamma}_j)}=\lambda(v_j)\bar{s}_j$,
and analogously: 
$s_\lambda\arrowvert_{U_\varrho(\bar{\Gamma}_j)+v}=\lambda(v_j)\bar{s}_j\circ\Phi(-v,\cdot)+\lambda(v)$.
Thus, $s_\lambda$ is a smooth function.

Choose an arbitrary representative $\omega'$ for
  $\lambda$. Since $\omega'$ is closed, the $1$-form
  $\tilde{p}^*\omega'\in \Omega^1(\tilde{M})$ is also closed, where
  $\tilde{p}:\tilde M\to M$ is the universal covering of $M$.
  Since each closed $1$-form on $\tilde M$ is exact, there exists
  $\tilde{f}\in C^\infty(\tilde{M})$ such that   $\tilde{p}^*\omega'=d\tilde{f}$. One can
  show easily that $\tilde{f}$ is invariant under the action of $h^{-1}(T)$ on
  $\tilde{M}$ and descends to $\bar f\in C^\infty(\bar{M})$,
  i.e., $\tilde{f}=\bar f\circ
  q$ where $q:\tilde{M}\to\tilde{M}/h^{-1}(T)=\bar{M}$ is the projection. We
  have $p\circ q=\tilde{p}$ and $q^*d\bar f=d\tilde{f}=\tilde{p}^*\omega'=
  q^*(p^*\omega')$ and hence $d\bar f=p^*\omega'$. Let
  $\bar{g}:=s_\lambda-\bar f\arrowvert_{U_\varrho(L)}:U_\varrho(L)\to\mathbb{R}$. A computation shows that for all $x \in U_\varrho(L)$ and $v
  \in H_1(M;\mathbb{Z})_{\mathbb{R}}$, we have $\bar{g}(x+v)=\bar{g}(x)$
and the existence of $g: U_\varrho(\Gamma)\to \mathbb{R}$ with $\bar{g}=g\circ p$ follows.

The map $g$ is smooth and we have on $U_\varrho(L)$:
\begin{align*}
p^*dg=d\bar{g}
=\sum_{i=1}^N\lambda(v_i)d\bar{s}_i-d\bar{f}
=p^*\left(\sum_{i=1}^N\lambda(v_i)ds_i-\omega'\right).
\end{align*}
Since $p$ is a surjective local diffeomorphism, the equality
$dg=\sum_{i=1}^N\lambda(v_i)ds_i-\omega'$ follows.

 Define now the smooth $1$-form
 \[\omega:=gd\zeta+(1-\zeta)\omega'+\zeta\sum_{i=1}^N\lambda(v_i)ds_i\]
 with $\zeta$ as in \eqref{stern}.
Using the fact that $\omega'$ is closed on $U_\varrho(\Gamma)$ and the
properties of $\zeta$, one can easily
verify that $\omega$ is smooth and closed.
Furthermore, for $x \in U_\varepsilon(\Gamma_j)$:
\begin{align*}
\omega\arrowvert_{x}=g(x)d\zeta\arrowvert_{x}+\left(1-\zeta(x)\right)\omega'\arrowvert_{x}+\zeta(x)\cdot\sum_{i=1}^N\lambda(v_i)ds_i\arrowvert_{x}
               =\lambda(v_j)ds_j\arrowvert_{x}, 
\end{align*}as claimed. We get
\begin{align*}
[\omega](v_j)=\int_{\gamma_j}\omega
=\lambda(v_j)\int_0^1ds_j\arrowvert_{\gamma_j(t)}\left(\dot{\gamma_j}(t)\right)dt\overset{(\ref{ds})}=\lambda(v_j)
\end{align*}
for $j=1,\dots,N$. With $\operatorname{span}\{v_1,\dots,v_N\}=H_1(M;\mathbb{R})$, this yields that
$\omega$ is a representative for $\lambda$.
\end{proof}
In the following, such a representative $\omega$ will be called a \emph{good
  representative} of $\lambda$ with respect to the family $\{v_1,\ldots,v_N\}$.

\section{Hedlund metrics}
Let $P$ be a centrally symmetric and convex polytope in $H_1(M;\mathbb{R})$ with nonempty interior, such that the directions of its vertices
are rational. Such a polytope will be called \emph{admissible}.  We call $\tilde{V}_P=\{\tilde{v}_1,\dots, \tilde{v}_N,-\tilde{v}_1,\dots,
-\tilde{v}_N\}$ the set of vertices  of $P$.

Let
$v_1,\dots,v_N$ be indivisible integer classes such that $v_i=\varepsilon_i\tilde{v}_i$ with $\varepsilon_i>0$,
$i=1,\dots,N$. 
Define $V_P:=\{v_1,\dots,v_N,-v_1,\dots,-v_N\}$ and let
$J_i$ be the subset of $V_P$ consisting of the indivisible integer classes corresponding to the
vertices belonging to the $i$-th
face $S_i$ of $P$. In order to simplify the notation, we assume without loss
of generality that
$J_1=\{v_1,\ldots,v_k\}$ for an integer $k\leq N$. 
The norm $|\cdot|$ on $H_1(M;\mathbb{R})$, whose unit ball is $P$, is given as follows
(for vectors lying in the cone over the face $S_1$):
\begin{equation}
v=\sum_{j=1}^{k}\alpha_j\tilde{v}_j \text{ with }\sum_{j=1}^{k}\alpha_j=1
\text{ and all }\alpha_j\geq0 \Rightarrow
|v|=1
\end{equation}
or generally
\begin{equation*}
v=\sum_{j=1}^{k}\alpha_j\tilde{v}_j \text{ with all }\alpha_j\geq0 \Rightarrow
|v|=\sum_{j=1}^{k}\alpha_j
\end{equation*}
and likewise for every other face of $P$.

Since $P$ is convex, for each face $S_i$ of $P$  exists  an element $\lambda_i$ of $H^1_{dR}(M)\simeq H^1(M,\mathbb{R})$
such that
\begin{align*}
\lambda_i(\tilde{v}_j)\left\{
\begin{array}{lll}
=&1,&\quad v_j=\varepsilon_j\tilde{v}_j\in J_i\\
<&1,&\quad v_j=\varepsilon_j\tilde{v}_j\not\in J_i\\
\end{array}\right.
\end{align*}
(i.e. $\lambda_i\equiv 1$ on the plane defined by the face $S_i$ and
$\lambda_i$ is
smaller on the rest of the polytope).
Now, since $P$ is symmetric, $-\lambda_i$ is the $1$-form corresponding to
$-S_i$ and we get in fact:
\begin{align}\label{lambda}
-1<\lambda_i(\tilde{v}_j)<1\text{ for }\pm v_j\not\in J_i.
\end{align}
We get an alternative definition for the norm:
\begin{equation}\label{defnorm}
v\in \bigoplus_{j=1}^k\mathbb{R}_{\geq 0}\cdot v_j\Rightarrow |v|=\lambda_1(v),
\end{equation}
and likewise for every other face of $P$.

The metrics defined below will be called \emph{Hedlund metrics} since such a
metric first appears in Hedlund's paper \cite{hedlund} in the case $M=\mathbb{T}^3$:

\begin{definition}\label{def}
Let  $P$ admissible polytope with vertices $\{\tilde
v_1,\ldots,\tilde v_N,-\tilde v_1,\ldots,-\tilde v_N\}$. 
Let $v_1,\ldots,v_N\in H_1(M,\mathbb{Z})_{\mathbb{R}}$ be
the indivisible integer 
classes such that
$\varepsilon_i\tilde{v}_i=v_i$ for some $\varepsilon_i>0$, $i=1,\dots,N$.
Choose disjoint admissible
curves $\gamma_1,\dots,\gamma_N$ representing the
classes  $v_1,\dots,v_N$.
For each face $S_i$ of $P$, let $\eta_i$ be a good representative of
$\lambda_i$ with respect to the family $\{v_1,\ldots,v_N\}$. 
A \emph{Hedlund metric} associated to $P$ on $(M,\rho)$ is a Riemannian metric $g$ that is conformal to $\rho$ and such that its dual metric $g^*$ satisfies: 
\begin{enumerate}
\renewcommand{\theenumi}{\arabic{enumi}}
\renewcommand{\labelenumi}{(H\,\theenumi)}
\item 
 $g^*_{\gamma_i(t)}(ds_i\arrowvert_{\gamma_i(t)},ds_i\arrowvert_{\gamma_i(t)})
=\underset{x\in U_\varepsilon(\Gamma_i)}\max\,
g^*_x(ds_i\arrowvert_x,ds_i\arrowvert_x)=\frac{1}{\varepsilon_i^2}$ for all
$t\in [0,1]$\newline
and
$g_x^*(ds_i\arrowvert_x,ds_i\arrowvert_x)<\frac{1}{\varepsilon_i^2}$
for $x \in U_\varepsilon(\Gamma_i)\setminus\Gamma_i$ and all $i \in \{1,\dots,N\}$.
\item $g^*_x(\eta_i\arrowvert_x,\eta_i\arrowvert_x)\leq 1$ for all $i=1,\ldots,N$ and $x\not\in U_\varepsilon(\Gamma)$.
\end{enumerate}
\end{definition}
Remark that for orientable compact surfaces of positive genus, it is not possible to choose disjoint loops
representing the vertices of the polytope. In fact, it is shown in \citet{ba}
that in the case of the $2$-torus, the stable norm induced by a Riemannian
metric on $\mathbb{T}^2$ has always a strictly convex unit ball. 
Yet, Massart shows in \cite{massart97} that this is not true in general: 
the stable norm induced by a smooth Finsler metric on a closed, orientable
surface has neither to be strictly convex, nor smooth. 
For a
\emph{non-orientable} surface, the analogon to Theorem \ref{main} can be found in
\citet{BaMa2008}: they show that if $M$ is a closed non-orientable
surface equipped with a Riemannian metric, then there exists in
every conformal class  a metric on $M$ whose stable norm has a polyhedron
as its unit ball.

\paragraph{Existence and properties of such a metric}
\begin{proposition}\label{propex}
On every compact Riemannian manifold $(M,\rho)$ with $\dim M\geq 3$ and for
every admissible $P$ in $H_1(M,\mathbb{R})$
there exists a Hedlund metric associated to $P$ on $(M,\rho)$.
\end{proposition}

\begin{proof}
Given the admissible polytope $P$, choose disjoint admissible
curves $\gamma_1,\dots,\gamma_N$ representing the indivisible integer classes
$v_1,\dots,v_N$ corresponding to its vertices \linebreak
$\tilde{v}_1,\dots,\tilde{v}_N$.
Let $\varepsilon_1,\dots,\varepsilon_N$ be the coefficients as in Definition
\ref{def}. 
For each face $S_i$ of $P$, $i=1,\dots,l$, let $\eta_i$ be a good representative for $\lambda_i$. 
Set \[\Omega:=\underset{\substack{j=1,\dots,l\\
                                 x \in M\setminus U_\varepsilon(\Gamma)}}\max
                             \rho^*_x(\eta_j\arrowvert_{x},\eta_j\arrowvert_{x})\] and 
 \[\Omega_i:=\max\{\underset{\substack{j=1,\dots,l\\x \in U_{\varrho}(\Gamma_i)}}\max
                             \frac{\rho^*_x(\eta_j\arrowvert_{x},\eta_j\arrowvert_{x})}{\rho^*_x(ds_i\arrowvert_{x},ds_i\arrowvert_{x})},\varepsilon_i^2\}\] for $i=1,\dots,N$. Define:
\begin{align*}
h_i: U_{\varrho}(\Gamma_i)&\to(0,\infty)\\
    x&\mapsto
    \frac{1}{\varepsilon_i^2\rho^*_x(ds_i\arrowvert_{x},ds_i\arrowvert_{x})}\cdot\exp(-C_i\cdot\ell(x)^2)
\end{align*}
 where
 \[
 C_i:=\ln\left(\frac{\Omega_i}{\varepsilon_i^2}\right)\cdot\frac{1}{\varepsilon^2}>0\]
and $\ell(x)$ is the distance from to $x$ to its ``projection''
$\gamma_i(s_i(x))\in \Gamma_i$.
Define the smooth function  $F:M\to(0,\infty)$ by
\begin{eqnarray*}
F(x)= \zeta(x)\cdot\sum_{i=1}^bh_i(x)+(1-\zeta(x))\cdot\frac{1}{\Omega},
\end{eqnarray*}
where $\zeta$ is a smooth bump function as in \eqref{stern}.  It is then easy to verify that the metric
$g$ defined by
\begin{equation*}
g^*_x=F(x)\rho^*_x  \quad \text{ for all } x \in M
\end{equation*}
is a Hedlund metric associated to $P$.
\end{proof}

\begin{proposition}
It results  immediately from Definition \ref{def} and from the properties of
an admissible polytope that 
\begin{equation}\label{normeta}
\|\eta_i\|^*:= \underset{x\in M}\max\big\|\eta_i\arrowvert_{x}\big\|^*_x=1
\end{equation}
for each face $S_i$ of $P$.
\end{proposition}

\begin{proof} Here again, we assume that $i=1$. The arguments are the same for
  every other face of $P$. 
Outside of  $U_\varepsilon(\Gamma)$, Definition \ref{def} yields
$\big\|\eta_1\arrowvert_{x}\big\|^*_{x}\leq 1$. With
\begin{eqnarray*}
\big\|\eta_1\arrowvert_{x}\big\|^*_{x}=\left\{
\begin{array}{ll}
\varepsilon_j\big\|ds_j\arrowvert_{x}\big\|^*_{x}=1 & ,x\in \Gamma_j \\
&\text{ and }j=1,\dots k\\
\varepsilon_j\big\|ds_j\arrowvert_{x}\big\|^*_{x}< 1& ,x\in
U_\varepsilon(\Gamma_j)\setminus\Gamma_j\\
&\text{ and }j=1,\dots k\\
\left|\lambda_1(v_j)\right|\cdot\big\|ds_j\arrowvert_{x}\big\|^*_x=  \varepsilon_j\left|\lambda_1(\tilde{v}_j)\right|\cdot\frac{1}{\varepsilon_j}\overset{(\ref{lambda})}<1\, &,x\in
U_\varepsilon(\Gamma_{j})\\
 &\text{ and }j>k,\\
\end{array}\right.
\end{eqnarray*}
this proves the statement.
\end{proof}

For the proof of the following lemma, we need to compute the lengths of the chosen
admissible curve $\gamma_1,\dots,\gamma_N$ relative to the new metric. Choose
$x=\gamma_i(t)\in\Gamma_i$ and a semi-geodesic chart $\varphi$ around
$x$. Recall the construction of such a chart; the matrix representing $\rho$
relative to the orthogonal basis
$(\dot{\gamma}_i(t),\partial_2^\varphi\arrowvert_{x},\dots,\partial_m^\varphi\arrowvert_{x})$
of $T_xM$ is diagonal. Hence, because $g$ is conformal to $\rho$, the
matrix representing $g$ relative to this basis is diagonal, too. Since the covectors
$(ds_i\arrowvert_{x},d\varphi_2\arrowvert_{x},\dots,d\varphi_m\arrowvert_{x})$
form a dual basis of $T_x^*M$, we obtain
\[g_x(\dot{\gamma}_i(t),\dot{\gamma}_i(t))=\frac{1}{g^*_x(ds_i\arrowvert_{x},ds_i\arrowvert_{x})},\]
using the fact that the matrice representing $g_x$ in the basis
$(\dot{\gamma}_i(t),\partial_2^\varphi\arrowvert_{x},\dots,\partial_m^\varphi\arrowvert_{x})$
is inverse to the matrix representing $g^*_x$ in the dual basis.
But because of (H $1$) in Definition \ref{def}, we have
$g^*_x(ds_i\arrowvert_{x},ds_i\arrowvert_{x})=\frac{1}{\varepsilon_i^2}$.
Hence, this leads to:
\begin{align}
L(\gamma_i)=\int_0^1\varepsilon_{i}dt=\varepsilon_{i}.\label{length}
\end{align}
It is possible to show that $\gamma_i$ is even the shortest curve representing
$v_i$: Assume, without loss of generality, that $v_i\in J_1$ and choose an
arbitrary curve $c:[0,1]\to M$ representing $v_i$. We have
$\lambda_1(v_i)=\varepsilon_i$ and hence
\begin{align*}
\varepsilon_i=&\int_c\eta_1=\int_0^1\eta_1\arrowvert_{c(t)}(\dot{c}(t))dt\leq\int_0^1\|\eta_1\arrowvert_{c(t)}\|^*\|\dot{c}(t)\|dt\\
&\overset{\eqref{normeta}}\leq\int_0^11\cdot\|\dot{c}(t)\|dt=L(c).
\end{align*} 

\begin{lemma}\label{obenab}
There is a constant $C=C(M,P)$ such that for each face $S_i$ of $P$, every $w\in\bigoplus_{v\in J_i}
\mathbb{N}\cdot v$
and every $x\in \bar{M}$, the distance from $x$ to $x+w$ is bounded abowe by $\lambda_i(w)+C$.
\end{lemma}

\begin{proof}
Recall the definitions of $\gamma_i$, $\Gamma_i$, $\bar{\gamma}_i$, $\bar{\Gamma}_i$, $i=1,\dots,N$, $L$ and $F_0$.
Define
\[D:=\underset{1\leq i,j\leq
  N}\max\,\,\underset{\substack{x\in\Gamma_i\\y\in\Gamma_j}}\min\,
d(x,y),\]

\[\operatorname{diam}(M):=\underset{x,y\in M}\max d(x,y)\]
and choose a real positive number $e$ such that $e>\underset{i=1,\dots,N}\max\varepsilon_i$.
Let
\begin{equation}\label{constant}
C:=2\cdot\operatorname{diam}(M) +\kappa\cdot (D+e)
\end{equation} 
where $d$ is the distance on $M$ induced from the Hedlund metric $g$ and 
$\kappa=\kappa(P)$ is the maximal number of vertices lying on a common face of
$P$.

Without loss of generality, we assume that $w\in\bigoplus_{v\in J_1}
\mathbb{N}\cdot v$, i.e. we can write $w=\sum_{i=1}^kn_iv_i$ with
$n_1,\dots,n_k\in\mathbb{N}$. 
We give a path from $x$ to $x+w$ that has length
bounded above by $\lambda_1(w)+C=\sum_{i=1}^k\varepsilon_in_i+C$. Assume that
$x\in F_0$ (otherwise, if $x\in F_0+u$ with $u\in H_1(M;\mathbb{Z})_{\mathbb{R}}$, we can replace the
path with startpoint $x-u$ as constructed below  with its image under $\Phi_u$). We join $x$ with $x+w$
by a path that runs \emph{as much as possible} in $L$ with ``changes of lines''
that are \emph{as short as possible}:

Choose $i_1\in\{j\mid 1\leq j\leq k,\, n_j\neq 0\}$ such that the point $x_1$ in $L\cap F_0$ with minimal
distance from $x$ lies in $\bar\Gamma_{i_1}$. Let $\tau_1$ be the corresponding geodesic segment from
$x$ to $x_1$ with minimal length. This length $\bar{L}(\gamma_1)$ is smaller
than $\operatorname{diam}(M)$. Let $c_1$ be the segment of
$\bar{\gamma}_{i_1}$ connecting $x_1$ and $x_1+n_{i_1}v_{i_1}$. This segment
  has length equal to 
\begin{align}
\bar L(c_1)=n_{i_1}\cdot  L(\gamma_{i_1})\overset{\eqref{length}}=n_{i_1}\cdot\varepsilon_{i_1}.\nonumber
\end{align}
Now choose $i_2\in \{j\mid 1\leq j\leq k,\, n_j\neq 0\}\setminus{i_1}$ and $x_2\in
\bar{\Gamma}_{i_2}+n_{i_1}v_{i_1}$ such that $x_2$ is the point of
$(L\setminus \bar{\Gamma}_{i_1})\cap
(F_0+n_{i_1}v_{i_1})$ having minimal distance from $\bar{\Gamma}_{i_1}\cap
(F_0+n_{i_1}v_{i_1})$. Let $x_1'$ be the point in $\bar{\Gamma}_{i_1}\cap
(F_0+n_{i_1}v_{i_1})$ at this minimal distance from $x_2$.
Let $c_1'$ be the section of  $\bar{\gamma}_{i_1}$ connecting $x_1$ and
$x_1'$; the length of $c_1'$ lies in
$[n_{i_1}\cdot\varepsilon_{i_1}-e,n_{i_1}\cdot\varepsilon_{i_1}+e]$. Let
$\tau_2$ be the minimal geodesic segment joining $x_1'$ and $x_2$, it has
length smaller than $D$. Now
continue in this way; choose $i_3\in \{j\mid 1\leq j\leq k,\, n_j\neq 0\}\setminus\{i_1,i_2\}$ and $x_3\in
\bar{\Gamma}_{i_3}+n_{i_1}v_{i_1}+n_{i_2}v_{i_2}$ such that $x_3$ is the point
of $(L\setminus(\bar{\Gamma}_{i_1}\cup\bar{\Gamma}_{i_2})) \cap
(F_0+n_{i_1}v_{i_1}+n_{i_2}v_{i_2})$ having minimal distance from $\bar{\Gamma}_{i_2}\cap
(F_0+n_{i_1}v_{i_1}+n_{i_2}v_{i_2})$. Let $x_2'$ be the point in $\bar{\Gamma}_{i_2}\cap
(F_0+n_{i_1}v_{i_1}+n_{i_2}v_{i_2})$ at this minimal distance from $x_3$. The curve $c_2'$
joining $x_2$ and $x_2'$ on $\bar{\Gamma}_{i_2}+n_{i_1}v_{i_1}$ has length
smaller than $n_{i_2}\cdot\varepsilon_{i_2}+e$.
 
If $n_j\neq 0$ for $j=1,\dots,k$, our path will be the composition 
\[\gamma:=\tau_1*c_1'*\tau_2*c_2'*\dots*c_k'*\tau_{i_{k+1}}\]
where $\tau_{k+1}$ is the path joining the last point in
$L\cap(F_0+\sum_{i=1}^kn_iv_i)$ with minimal distance from $x+w$ to $x+w$ and
has length smaller than $\operatorname{diam}(M)$. Summing all
the lengths of those segments we get 
\begin{align*}
\bar L(\gamma)\leq&\operatorname{diam}(M)+n_{i_1}\cdot\varepsilon_{i_1}+e+D+n_{i_2}\cdot\varepsilon_{i_2}+e+D\\&+\dots+n_{i_k}\cdot\varepsilon_{i_k}+e+\operatorname{diam}(M)\\
&=\lambda_1(w)+k\cdot
e+k\cdot D+2\cdot \operatorname{diam}(M)\leq \lambda_1(w)+C.
\end{align*}
Finally, if $n_j=0$ for some $j\in\{1,\dots,k\}$, we need to make fewer
changes of lines, and the inequality can be shown the same way.
\end{proof}

\paragraph{The stable norm and the main theorem.}
In the introduction of this paper, we gave the definition of the stable
norm induced from a Riemannian metric $g$ on $M$. Here we give a way to compute
 the stable norm of a vector lying in $H_1(M;\mathbb{Z})_{\mathbb{R}}$:
Define \begin{eqnarray*}
f:H_1(M;\mathbb{Z})_{\mathbb{R}}&\to&\mathbb{R}_{\geq0}\\
v&\mapsto&\inf\{L(\gamma)|\gamma \text{ closed curve representing }
v\}
\end{eqnarray*}
and $f_n:n^{-1}H_1(M;\mathbb{Z})_{\mathbb{R}}\to
\mathbb{R}_{\geq0}$, $f_n(v)=n^{-1}f(nv)$.
In \citet{ba} it  is shown that $f_n$ converges uniformly on compact sets to the
stable norm $\|\cdot\|_s$. Especially, we have: if $(v_n)_{n \in \mathbb{N}}$ is a sequence in
$H_1(M;\mathbb{Z})_{\mathbb{R}}$ with
$\lim_{n\to\infty}\frac{v_n}{n}=v \in H_1(M;\mathbb{R})$ (relative to the
standard topology on the vector space $H_1(M;\mathbb{R})\simeq\mathbb{R}^b$), then we have for
the norm of $v$:
\begin{equation*}
\|v\|_s=\underset{n\to\infty}\lim \frac{f(v_n)}{n}.
\end{equation*}
If $\bar{d}$ is the distance on $\bar{M}$ induced from $p^*g$,  we have
for $v \in H_1(M;\mathbb{Z})_{\mathbb{R}}$:
\begin{equation*}f(v)=\underset{x\in\bar{M}}\inf
\bar{d}(x,x+v)=\underset{x\in F_0}\min\, \bar{d}(x,x+v)
\end{equation*} because $p^*g$ is a periodic metric and the closure of $F_0$
is a compact set. 
With $\underset{n\to\infty}\lim \frac{nv}{n}=v$, this yields:
\begin{eqnarray*}
\|v\|_s&=&\underset{n\to\infty}\lim \frac{f(nv)}{n}\\
     &=&\underset{n\to\infty}\lim\frac{\underset{x\in F_0}\min\, \bar{d}(x,x+nv)}{n}.
\end{eqnarray*}

\begin{theorem}\label{main}
The polytope $P$ is the unit ball of the stable norm on $H_1(M;\mathbb{R})$
induced by  an arbitrary Hedlund metric associated to $P$ on $M$.
\end{theorem}

Note that by Definition \ref{def}, the Hedlund metric is chosen in
the conformal class of  the given Riemannian
metric $\rho$ on $M$.
\begin{proof}
Let $g$ be a Hedlund-metric associated to $P$. We show that for each $w\in
\bigoplus_{j=1}^k\mathbb{N}\cdot v_j$, the stable norm of $w$ is given by
$\|w\|_s=\lambda_1(w)$. The proof of this works analogously for every other
face of $P$.  Consequently, this holds for all vectors in $H_1(M;\mathbb{R})$
that can be written as linear combinations of the vectors $v_1,\dots,v_N$ with
rational coefficients, and then, by continuity, this holds for all vectors in $H_1(M;\mathbb{R})$.
Let $x$ be an arbitrary point in $F_0$ and let $n\in \mathbb{N}$. Let
$\gamma:[0,1]\to \bar{M}$ be an arbitrary path from $x$ to $x+nw$. We have 
\begin{align*}
\lambda_1(nw)&=\int_\gamma\eta_1=\int_0^1\eta_1\arrowvert_{\gamma(t)}(\dot{\gamma}(t))dt\leq\int_0^1\|\eta_1\arrowvert_{\gamma(t)}\|^*\|\dot{\gamma}(t)\|dt\\
&\overset{\eqref{normeta}}\leq\int_0^11\cdot\|\dot{\gamma}(t)\|dt=\bar{L}(\gamma)
\end{align*}
With this and Lemma \ref{obenab} we get 
\begin{equation*}
\lambda_1(n\cdot w)\leq \bar{d}(x,x+nw)\leq \lambda_1(n\cdot w)+C.
\end{equation*} 
Thus 
\begin{equation*}
\lambda_1(n\cdot w)\leq \underset{x\in F_0}\min\,\bar{d}(x,x+nw)\leq \lambda_1(n\cdot w)+C,
\end{equation*}
and 
\begin{equation*}
\lambda_1(w)\leq \frac{\underset{x\in F_0}\min\,\bar{d}(x,x+nw)}{n}\leq \lambda_1(w)+\frac{C}{n}.
\end{equation*}
Letting $n$ go to infinity, this yields $\|w\|_s=\lambda_1(w)$, as claimed.
\end{proof}

\end{document}